\documentstyle{article}
\author{Ehud Hrushovski \footnote{Hebrew University at Jerusalem.  Supported
by the Miller Institute, University of California, Berkeley}
\\ Thomas Scanlon \footnote{Mathematical Sciences Research Institute.  Supported by an NSF MSPRF}}
\title{Lascar and Morley ranks differ in differentially closed fields}

\newtheorem{theorem}{Theorem}[section]

\newtheorem{remark}[theorem]{Remark}

\newcommand{\proof}{\noindent \it Proof \hbox{       }}

\font\msbmten msbm10
\def \Bbb#1{\hbox{\msbmten {#1} }}
\def \iso \simeq
\def\acl{{\rm acl}}

\def\k{\bf \it k}
\def\proof{ { {\noindent \it Proof} \hbox{       }  }  }

\def\iso{\simeq}
\def\acl{{\rm acl }}

\def\Uu{\Bbb U}

\def\tp{ \rm{tp}}
\def\<{\begin}
\def\>{\end}

\begin{document}
\maketitle

We note here, in answer to a question of Poizat,
that the Morley and Lascar ranks need not coincide in differentially
closed fields.  We will approach this through  the (perhaps) 
more fundamental issue of the variation of Morley rank in families.  
We will be interested here only in sets of finite Morley rank.
\S~\ref{DefMR}
consists of some general lemmas relating the above issues. 
\S~\ref{nondef}
points out a family of sets of finite Morley rank, whose
Morley rank exhibits  discontinuous upward jumps.  To make the
base of the family itself have finite Morley rank, we use a theorem 
of Buium.

We thank John Baldwin, Anand Pillay, and Wai Yan Pong for reading
an earlier version of this note and suggesting improvements.

\<{section} {Definability of Morley rank}
\label{DefMR}

We will say that Morley rank is \emph{definable} 
(respectively \emph{upward}, resp.
\emph{downward semi-definable}) if for every set of parameters $A$ and
 $A$-definable family of definable
sets $E_b (b \in B)$ and $b \in B$, there is an $A$-definable
set $B' \subseteq B$ such that $b \in B'$ and 
$MR(E_{b'})= MR(E_b) $  (resp. $\geq$, $\leq$) for $b' \in B'$.

\<{lemma} \label{defMR} Let $T$ be a theory of finite Morley rank.  If  
 Morley and Lascar ranks coincide on definable sets,
 then Morley rank is downward semi-definable.
 \>{lemma}

\proof   Suppose  Morley rank is not downward semi-definable,
 and let $E_b$ ($b \in B$)
be an $A$-definable  family demonstrating this. 
That is, there is some $b^* \in B$
so that for every $A$-definable set $B' \subseteq B$ with $b \in B'$
there is some $b' \in B'$ with $MR(E_{b'}) > MR(E_{b^*})$.  
Replace $B$ with an $A$-definable set of minimal Morley rank and 
degree containing $b$; so that now $MR(\tp (b^*/A)) =  MR(B) := m$.
Let $d := MR(E_{b^*})$.
For $b \in B$, let $E'_b := {E_b}^m$.  Then 
$MR(E'_b) = mMR(E_b)$.  So it is always a multiple of $m$.  We have
$MR(E'_{b^*})=md$,
while for many $b'$, $MR(E'_{b'}) \geq md+m$.

For $B' \subset B$, let 
$$X_{B'} := \{(e,b):  b \in B', e \in E'_b \}$$
and $X := X_B$.  Replace $B$ with some $A$-definable $B' \subseteq
B$ with $b^* \in B'$ and the property that for every  $A$-definable
$B'' \subseteq B'$ with $b^* \in B''$, $(MR(X_{B''}), dM(X_{B''}) )
= (MR(X_{B'}), dM(X_{B'}) )$.

If $(e,b) \in X$ and $b$ is not generic in
$B$, then $b \in B'$ for some $A$-definable subset of $B$ with 
$MR(B') < MR(B)$ and $(e,b) \in X_{B'}$. Since 
$X_{B \setminus B'}$ and $X_B$ have the same
Morley rank and degree by the above reduction,  
$MR(X_B') < MR(X)$; so, $MR(\tp(e,b/A)) < MR(X)$.  

On the other hand, if $(e,b) \in X$, and $b \in B$ is generic,
then $U(\tp(e,b/A)) \leq md+m$.
But $MR(E'_b) \geq md+m$ for infinitely many $b \in B$;
 so $MR(X) > md+m$.
Thus $U(\tp(e,b/A)) < MR(X)$ for any $(e,b) \in X$.
So $U(X) < MR(X)$, a contradiction.   

\>{section}

\<{section}{A non-definable family }
\label{nondef}

We now work with differential fields of characteristic $0$, and fix
a universal domain $\Uu$ (a saturated differentially closed field.)

Our plan is to produce a finite rank
definable family of abelian varieties whose
Manin kernels exhibit non-definable jumps in
Morley rank.  One difficulty is that there does
not exist a definable family of abelian varieties
containing a copy of every abelian variety of a given 
dimension.  However, there are definable families of
abelian varieties  containing isomorphic copies of 
every \emph{principally polarized} abelian variety of a given dimension.

For every abelian variety $A$ there 
is another abelian variety $\check{A}$, called the 
dual abelian variety, which parametrizes the line 
bundles on $A$ algebraically equivalent to zero.
A polarization is an isogeny $\lambda: A \rightarrow 
\check{A}$. A polarization is principal if it is an 
isomorphism. A \emph{principally polarized abelian variety}
is an abelian variety $A$ given together with a principal
polarization $\lambda: A \rightarrow \check{A}$. 
Not all abelian varieties admit a principal polarization,
but elliptic curves  
always do.

\<{theorem}[\cite{GIT} VII \S2]\label{universal}
Let $L$ be an algebraically closed
field, and $g$ a positive integer.
There exists a definable family 
$\{ (A_b, \lambda_b) : \ b \in F \}$ of $g$-dimensional 
principally polarized abelian varieties, such that every 
principally polarized $g$-dimensional abelian variety over $L$ is
isomorphic to some $(A_\alpha, \lambda_\alpha)$.   \>{theorem}

Let $E_t$ be an elliptic curve  with $j$-invariant $t$.
Given $t, t'$, let $E(t,t') := E_t \times E_{t'}$.  Given also
 an integer $n$,
there exist group-theoretic isomorphisms $\iota $ between
 the finite $n$-torsion 
subgroups of $E_t$ and of $E_{t'}$.
 Let $E(t,t', \iota, n)$ be the quotient
 of $E_t \times E_{t'}$ by
the graph of $\iota$.  Then $A :=  E(t,t',\iota, n)$ is an Abelian 
variety of dimension $2$.
When $E_t,E_{t'}$ are not isogenous, 
$A$ has precisely two connected definable subgroups of
Morley rank 1, namely the images of $E_t$ and of $E_{t'}$.
Their intersection has
order $n^2$.
For a general choice of $\iota$, $E(t,t', \iota, n)$ 
need not admit a principal polarization.
However, if we choose $\iota$ to be 
anti-symplectic (ie 
$\langle \iota(x) , \iota (y) \rangle_{E_{t'}} =
\langle y, x \rangle_{E_{t}}$ where $\langle \cdot, \cdot \rangle$ is
the Weil pairing (See~\cite{Mil} \S 16 for the 
general theory of the Weil pairing) then $E(t,t', \iota, n)$ is self-dual.

\<{lemma} \label{ppav}  If $\iota : E_t[n] \rightarrow E_{t'}[n]$
is an anti-symplectic isomorphism of the $n$-torsion points, then
$A := E(t,t',\iota, n)$ has a natural principal polarization. \>{lemma}
 
\proof  
Since $\iota$ is an anti-symplectic map, the graph of 
$\iota$ is isotropic for the pairing on $E_t \times 
E_{t'}$ :
\begin{eqnarray*}
\langle (x,\iota(x)), (y, \iota(y)) \rangle_{E_t \times 
E_{t'}} & = & \langle x, y \rangle_{E_t} \cdot \langle \iota(x),
		\iota(y) \rangle_{E_{t'}} \\
	& = & \langle x, y \rangle_{E_t} \langle y, x \rangle_{E_{t}} \\
	& = & 1 
\end{eqnarray*}

Since $\# (E_t \times E_{t'}) [n] = n^4$ and the pairing is perfect,
a maximal isotropic space has size $n^2$ which is the size of the 
graph of $\iota$.  Hence, the graph of $\iota$ is a maximal isotropic 
subspace.

The lemma now follows from the more general lemma:

\<{lemma} \label{ppiso} 
Let $A$ be a principally polarized abelian variety identified 
	with its dual via the polarization.  Let $\Gamma \subseteq
	A[n]$ be a maximal isotropic subgroup of the $n$-torsion 
	subgroup of $A$.  Let $B = A/\Gamma$.  Then $B$ also admits a 
	principal polarization. \>{lemma}

\proof  Let $\pi: A \to B$ be the quotient map.  Let 
	$\phi: B \to A$ be defined by $b \mapsto [n] a$ 
	where $a$ is any choice of a pre-image of $b$ under
	$\pi$.  $\phi$ induces a dual exact sequence
$$ 0 \longrightarrow \textrm{Hom} (\ker{\phi}, \mu_n)
 \longrightarrow A \longrightarrow \check{B} 
	\longrightarrow 0$$
where $\mu_n$ is the group of $n$-th roots of unity 
(See~\cite{Mum} III \S 15).  
The kernel of $\phi$ is $A[n]/\Gamma$.  Since $\Gamma$ is 
isotropic, the pairing on $A[n]$ descends to a pairing 
$(A[n]/\Gamma) \times \Gamma \rightarrow \mu_n$.  Since 
the pairing is perfect and $\# A[n] = (\# \Gamma)^2$, 
via this pairing $\Gamma = \textrm{Hom} (A[n]/\Gamma, \mu_n)$.

Thus, the above exact sequence is 
$$0 \longrightarrow \Gamma \longrightarrow A
 \longrightarrow \check{B} \longrightarrow 0 $$

That is, $\check{B}$ is the quotient $A/\Gamma = B$.

\<{lemma}\label{simp}
Let $F'$ be a Zariski (resp. Kolchin) closed subset of $F$, the 
definable parameter space for two dimensional principally polarized
 abelian varieties
of Theorem~\ref{universal}.
  Assume $F'$ has a Zariski (resp. Kolchin) dense
subset $\{t_1,t_2,\ldots\}$, such that $A_{t_n}$ is isomorphic 
to some
$E(t,t',\iota, n )$  with $\iota$ anti-symplectic as above.
 Then for generic $t \in F'$, $A_t$ is a simple
abelian variety.   \>{lemma}

\proof  Otherwise, a generic $A_t$ contains
two elliptic curves.  Their   intersection is necessarily finite,
 say of order
$m$.  But then infinitely many $A_{t_n}$ must contain two elliptic
 curves 
with intersection of order $m$.  For $n>m$, this contradicts 
the remarks above. 

At this point it is quite easy to see that exists in DCF$_0$
 a definable
family of definable sets, whose generic element is strongly 
minimal, but
with densely many sets of Morley rank 2.  Thus:

\<{cor}  In $DCF_0$, Morley rank is not downwards semi-definable.
 \>{cor}

However, since $DCF_0$ does not have finite Morley
rank, lemma \ref{defMR} does not directly apply.  
At this point we quote a theorem from \cite{Bui}.

\<{theorem}[Buium~\cite{Bui}]\label{buium}  Let $(A, \lambda)$ be
 any principally polarized
 abelian variety
 of maximal $\delta$-rank.
  There exists
a definable family $\{ (A_t, \lambda_t): \ t \in F_1 \}$,
 containing a definably isomorphic
copy of every principally polarized abelian variety isogenous to $A$,
and such that $F_1$
has finite Morley rank.   \>{theorem}

We leave the notion of \emph{$\delta$-rank} undefined here since 
we need only the facts that:
\begin{itemize}
\item A generic elliptic curve has maximal $\delta$-rank.
\item The property of having maximal $\delta$-rank is isogeny invariant.
\item The product of two abelian varieties each of maximal 
	$\delta$-rank is also of maximal $\delta$-rank.
\end{itemize}

It seems likely that the $\delta$-rank condition is unnecessary in 
Buium's theorem, but we leave this issue aside.
 
\<{cor}  There exists a finite Morley rank definable subset
$Y$, such that Morley rank is not 
downwards semi-definable. \>{cor}

\proof  Pick $t,t'$ algebraically independent over $\k$, the 
field of differential constants of $\Uu$.  Let
$J_t, J_{t'}$ be elliptic curves with $j$-invariants $t, t'$.  Let
$A := J_t \times J_{t'}$, and let $F_1$ be a family as guaranteed to
exist by Theorem~\ref{buium}.  Given $n$, pick $c=c(n) \in F_1$ with 
$A_c$ isomorphic to $E(t,t', \iota, n)$.  Let $F_2$ be the Kolchin
closure of the set $\{c(1),c(2),...\}$.  Let $b$ be a generic element
of $F_2$.  By Lemma~\ref{simp}, $A_{b}$ is a simple Abelian
variety.  If $A_b$ were isogenous to an Abelian variety defined over
$\k$, this would be guaranteed by a certain formula true of $b$,
and the same formula would hold of infinitely many $c(n)$; hence
$A$ would also have this property, contradicting the choice of $t, t'$.
Thus $A_b$ is a simple, non-isotrivial Abelian variety.

For $t \in F_2$, let  $M_t$ be the Manin kernel of $A_t$.
$M_t$ is uniformly definable over $t$ (cf. \cite{hru-ita}). 
 Then (cf. \cite{hru-sok}) $M_t$ has Morley rank 1 for generic 
$t \in F'$
(when $A_t$ is a nonisotrivial simple Abelian variety.)
  But it has Morley 
rank 2 for each $t=c(n)$
(when $A_t$ is isogenous to a product of elliptic curves.)
  Thus Morley rank is not downward semi-definable in 
$Y = \{(a,t) : t \in F^\#, a \in M_t \}$.

\<{cor} Morley and Lascar rank do not agree on definable sets
 in DCF$_0$. \>{cor}
\proof Since $Y$ has finite Morley rank, with the structure
 induced from 
the ambient differentially closed field,
 Lemma~\ref{defMR} applies.

\<{question} \ 
\>{question}
Marker and Pillay have noted that on $0$-definable sets of
 differential order $2$,
Lascar and Morley ranks are the same.
  Examples similar to the one produced 
above have order at least $5$.  Is there a theorem 
responsible for this gap?

\>{section}

\>{document}